\documentclass[10pt]{article}
\usepackage{amssymb, amsmath, url, graphicx, setspace, geometry}
\usepackage{calrsfs}
\usepackage{wasysym, xcolor}

\def\3{\subset }
\def\4{\subseteq }
\def\<{\left<}
\def\>{\right>}

\def\bit{\begin{itemize}}
\def\eit{\end{itemize}}
\def\3{\subset }
\def\4{\subseteq }

\def\0{\leqno}

\def\barr{\begin{array}}
\def\earr{\end{array}}

\def\Z{{\rlap{$\kern2pt{\rm Z}$}{\rm Z}\,}}
\def\bld#1#2{{\buildrel{#1}\over{#2}}}
\def\st#1#2{{\mathrel{\mathop{#2}\limits_{#1}}{}\!}}
\def\stb#1#2#3{{\st{{#1}}{\bld{{#2}}{#3}}{}\!}}
\def\xmare#1#2{\stb{#1}{#2}{\mbox{\Large$\times$}}}

\DeclareMathOperator{\cdeg}{cdeg}

\title{\bf Element orders in extraspecial groups}
\author{Mihai-Silviu Lazorec}
\date{April 14, 2024}

\begin{document}

\maketitle

\begin{abstract}
By using the structure and some properties of extraspecial and generalized/almost extraspecial $p$-groups, we explicitly determine the number of elements of specific orders in such groups. As a consequence, one may find the number of cyclic subgroups of any (generalized/almost) extraspecial group. For a finite group $G$, the ratio of the number of cyclic subgroups to the number of subgroups is called the cyclicity degree of $G$ and is denoted by $\cdeg(G)$. We show that the set containing the cyclicity degrees of all finite groups is dense in $[0, 1]$. This is equivalent to giving an affirmative answer to the following question posed by T\'{o}th and T\u arn\u auceanu: \textit{``For every $a\in [0, 1]$, does there exist a sequence $(G_n)_{n\geq 1}$ of finite groups such that $\displaystyle\lim_{n\to\infty} \cdeg(G_n)=a$?"}. We show that such sequences are formed of finite direct products of extraspecial groups of a specific type. 
\end{abstract}



\section{Introduction}

Let $G$ be a finite group and $C_n$ be the cyclic group with $n$ elements where $n\geq 1$ is an integer. We denote by $o(x)$, $L(G)$ and $C(G)$ the order of an element $x$ of $G$, the subgroup lattice of $G$ and the poset of cyclic subgroups of $G$, respectively. The research done on the element orders of $G$ has strong connections with other intensively studied problems of finite group theory. Some of these are:
\begin{itemize}
\item[$(P_1)$] characterizing the nature and structure of $G$ via various tools defined using element orders;
\item[$(P_2)$] recognizing $G$;
\item[$(P_3)$] evaluating the number of (cyclic) subgroups of $G$;
\item[$(P_4)$] estimating the number of solutions of an equation in $G$.
\end{itemize}

In what concerns $(P_1)$, maybe the most popular topic related to it is the study of the so-called sum of element orders of a finite group $G$. This sum is denoted by $\psi(G)$ and is defined as
$$\psi(G)=\sum\limits_{x\in G}o(x).$$
During the last two decades a lot of papers focused on the properties of $\psi(G)$. Just to indicate a few, one may want to check \cite{1, 3, 4, 32}; also, \cite{15} is a recent survey about this topic. 

Shi's conjecture (see \cite{29}) stating that each finite simple group $G$ is determined by its order $|G|$ and its spectrum $\omega(G)$ (i.e. its set of element orders) was validated in \cite{36}. This is a remarkable result related to $(P_2)$. It set a fruitful research topic and we recommend the reader to check \cite{12} for a recent survey.

In regard to $(P_3)$, it is well-known that an element $x$ of a finite group $G$ generates a cyclic subgroup $\langle x\rangle\cong C_{o(x)}$ of $G$. Also, there is a strong connection between the properties of the  subgroup lattice $L(G)$ of $G$ and $G$ itself. The monograph \cite{28} outlines a lot of results which support this idea. Also, evaluating the number of subgroups of $G$, especially when $G$ is a $p$-group, is a research topic which goes back to old papers such as \cite{9, 20}. It is  worth mentioning that this subject is still of great interest and one may consult more recent papers (see \cite{5, 10, 27, 33}).
  
Finally, Frobenius's theorem (see \cite{9}) stating that if $n$ divides the order $|G|$ of a finite group $G$, then the number of solutions of the equation $x^n=1$ in $G$ is a multiple of $n$, probably is the most notable result associated with $(P_4)$. We refer the reader to \cite{14, 18, 19} for various proofs, applications and generalizations of this result.

We recall that a finite $p$-group $G$ is called:
\begin{itemize}
\item[--] special if either $G$ is elementary abelian or $G$ is of class 2 and $G'=Z(G)=\Phi(G)$ is elementary abelian (p.183 of \cite{11});
\item[--] extraspecial if $G$ is special and $|G'|=p$ (p.183 of \cite{11});
\item[--] almost extraspecial if $G'=\Phi(G)$, $|G'|=p$ and $Z(G)\cong C_{p^2}$ (Definition 2.2 of \cite{7});
\item[--] generalized extraspecial if $G'=\Phi(G)$, $|G'|=p$ and $G'\leq Z(G)$ (Definition 3.1 of \cite{30}).
\end{itemize}
Also, for a finite group $G$, the authors of \cite{34} define the cyclicity degree of $G$, denoted by $\cdeg(G)$, as being the quantity
$$\cdeg(G)=\frac{|C(G)|}{|L(G)|}.$$
This ratio measures the probability that randomly selecting a subgroup $H$ of $G$, $H$ is cyclic. In the same paper, the authors determine various properties of the cyclicity degree and explicit formulas for $\cdeg(G)$ where $G$ belongs to some classes of finite groups. They also pose the following question:\\ 

\textbf{Question 1.1.} \textit{For every $a\in [0, 1]$, does there exist a sequence $(G_n)_{n\geq 1}$ of finite groups such that $\displaystyle\lim_{n\to\infty} \cdeg(G_n)=a$?}\\

We conclude this section by describing the organization and purposes of this manuscript. The origin of this paper is based on some of the author's previous work. More exactly Question 2 of \cite{21} asks for an explicit formula of $\psi(G)$ if $G$ is a non-abelian special 2-group. This problem remains open, but can be solved for extraspecial $p$-groups, in particular. We also determine an explicit formula of $\psi(G)$ if $G$ is a generalized/almost extraspecial $p$-group. These results are obtained by counting the number of elements of specific orders in a (generalized/almost) extraspecial $p$-group. Section 2 of the paper deals with this first objective. Our second goal is to provide an affirmative answer to Question 1.1. This is done in Section 3, where we prove that each term of the sequence $(G_n)_{n\geq 1}$ can be formed by using extraspecial groups of a specific type.    
  
\section{On the element orders in (generalized/almost) extraspecial \textit{p}-groups}

Let $p$ be a prime, $G$ and $H$ be finite $p$-groups and $k\geq 0, n\geq 2$ be integers. We denote the central product of $G$ and $H$ by $G\circ H$. The central product of $n$ copies of $G$ is denoted by $G^{\circ n}$. By convention, one has $G^{\circ 0}=\lbrace 1\rbrace$ and $G^{\circ 1}=G$. The exponent of $G$ is denoted by $\exp(G)$, while $n_{p^k}(G)$ is the number of elements of order $p^k$ of $G$. We also consider the sets
$$\Omega_{\lbrace k\rbrace}(G)=\lbrace x\in G \ | \ o(x)|p^k\rbrace$$ which generate the omega subgroups $\Omega_k(G)$ of $G$ while $$\mho^k(G)=\langle\lbrace g^{p^k} \ | \ g\in G\rbrace\rangle$$ are the agemo subgroups of $G$. The dihedral group with 8 elements and the quaternion group are denoted by $D_8$ and $Q_8$, respectively. The modular $p$-group of order $p^n$ is
$$M_{p^n}=\langle x, y \ | \ x^{p^{n-1}}=y^p=1, y^{-1}xy=x^{p^{n-2}+1}\rangle,$$
where $n\geq 3$ if $p$ is odd, and $n\geq 4$ if $p=2$. Finally, for an odd $p$, the Heisenberg group modulo $p$ is of order $p^3$ and its representation is
$$He_p=\langle x, y, z \ | \ x^p=y^p=z^p=1, [x, z]=[y, z]=1, [x, y]=z\rangle.$$

Central products play an important role in understanding the structure of a (generalized/almost) extraspecial $p$-group. By summarizing the results given in section 3.3 (especially Theorem 3.14) of \cite{8}, section 2 (especially Theorem 2.3) of \cite{7} and Lemma 3.2 of \cite{30}, one gets a complete classification of the (generalized/almost) extraspecial $p$-groups.\\

\textbf{Lemma 2.1.} \textit{Let $G$ be a finite $p$-group.
\begin{itemize}
\item[i)] If $G$ is extraspecial, then there is an integer $n\geq 1$ such that $|G|=p^{2n+1}$ and \\
$\bullet$ $G\cong D_8^{\circ n}$ or $G\cong D_8^{\circ (n-1)}\circ Q_8$ if $p=2$;\\
$\bullet$ $G\cong He_p^{\circ n}$ or $G\cong M_{p^3}^{\circ n}$ if $p$ is odd.
\item[ii)] If $G$ is almost extraspecial, then there is an integer $n\geq 1$ such that $|G|=p^{2n+2}$ and \\
$\bullet$ $G\cong D_8^{\circ n}\circ C_4$ if $p=2$;\\
$\bullet$ $G\cong He_p^{\circ n}\circ C_{p^2}$ if $p$ is odd.
\item[iii)] If $G$ is generalized extraspecial, then $G\cong E\times A$ or $G\cong (E\circ C_{p^2})\times A$ where $E$ is extraspecial and $A$ is elementary abelian; more exactly, there are some integers $n\geq 1, k\geq 0$ and\\
$\bullet$ $G\cong D_8^{\circ n}\times C_2^k$ or $G\cong (D_8^{\circ (n-1)}\circ Q_8)\times C_2^k$ or $G\cong (D_8^{\circ n}\circ C_4)\times C_2^k$ if $p=2$;\\
$\bullet$ $G\cong He_p^{\circ n}\times C_p^k$ or $G\cong M_{p^3}^{\circ n}\times C_p^k$ or $G\cong (He_p^{\circ n}\circ C_{p^2})\times C_p^k$ if $p$ is odd.
\end{itemize}}

In what concerns item \textit{iii)} of Lemma 2.1, note that even though there are 4 isomorphism types of extraspecial $p$-groups, there are only 3 isomorphism types of generalized extraspecial $p$-groups. This is a consequence of the facts that 
$$D_8\circ C_4\cong Q_8\circ C_4 \text{ \ and \ } He_p\circ C_{p^2}\cong M_{p^3}\circ C_{p^2}$$ 
since these are leading to 
$$D_8^{\circ n}\circ C_4\cong (D_8^{\circ (n-1)}\circ Q_8)\circ C_4 \text{ \ and \ } He_p^{\circ n}\circ C_{p^2}\cong M_{p^3}^{\circ n}\circ C_{p^2},$$
respectively. 

If one is interested in studying the element orders of a finite $p$-group $G$, then it is helpful to know preliminary information about $\exp(G)$. Let $G$ be a (generalized/almost) extraspecial $p$-group. Since $|\Phi(G)|=p$, it follows that $\exp(G)\in\lbrace p, p^2\rbrace$. By Theorem 5.2 (ii) of \cite{11}, it is known that if $G\cong He_p^{\circ n}$, then $\exp(G)=p$. It follows that $\exp(G)=p$ if $G\cong He_p^{\circ n}\times C_p^k$, as well. If $G$ is isomorphic to any other group outlined in Lemma 2.1, then by the construction of the central product, it is known that $G$ has a subgroup $H$ such that $H\cong D_8$ or $H\cong Q_8$ if $p=2$, and $H\cong M_{p^3}$ or $H\cong C_{p^2}$ if $p$ is odd. In either case, we have that $\exp(H)=p^2$ and, since $\exp(H)|\exp(G)$ and $\exp(G)\leq p^2$, it follows that $\exp(G)=p^2$. To summarize, we outline another preliminary result. \\

\textbf{Lemma 2.2.} \textit{Let $G$ be a (generalized/almost) extraspecial $p$-group and $n\geq 1, k\geq 0$ be integers.
\begin{itemize}
\item[i)] If $p=2$, then $\exp(G)=4$;
\item[ii)] If $p$ is odd, then
$\exp(G)=\begin{cases} p &\mbox{, if } G\cong  He_p^{\circ n} \text{ or } G\cong He_p^{\circ n}\times C_p^k\\ p^2 &\mbox{, in any other cases}\end{cases}.$
\end{itemize}}
 
Since a generalized extraspecial $p$-group $G$ is a direct product of finite $p$-groups, it is useful to count the number of elements of order $p$ of a direct product with respect to the numbers of elements of order $p$ of its components.\\ 

\textbf{Proposition 2.3.} \textit{Let $G\cong G_1\times G_2\times\ldots\times G_l$ where $l\geq 2$ is an integer and $G_i$ is a finite group for all $i\in\lbrace 1,2,\ldots, l\rbrace$. Then 
$$n_p(G)=\prod\limits_{i=1}^l (n_p(G_i)+1)-1.$$}

\textbf{Proof.} Under the above hypotheses, let $x=(x_1, x_2,\ldots ,x_l)\in G$. Then 
$$o(x)=p \Longleftrightarrow \begin{cases} x^p=1\\ x\not =1\end{cases}\Longleftrightarrow\begin{cases} x_i^p=1_{G_i}, \ \forall \ i\in\lbrace 1,2,\ldots, l\rbrace\\ \exists \ i\in\lbrace 1,2,\ldots, l\rbrace \text{ such that } x_i\not =1_{G_i}\end{cases},$$
where $1_{G_i}$ denotes the identity of $G_i$ for all $i\in\lbrace 1,2,\ldots, l\rbrace$. Hence, the general form of a solution of the last system is $x=(x_1, x_2, \ldots, x_l)$ where at least one component $x_i$ is an element of order $p$ in $G_i$. Therefore, the number of solutions (which coincides with $n_p(G)$) is $$n_p(G)=\prod\limits_{i=1}^l (n_p(G_i)+1)-1,$$
as desired.
\hfill\rule{1,5mm}{1,5mm}\\

In order to study the element orders in a (generalized/almost) extraspecial $p$-group, we distinguish two cases according to the parity of $p$. We first suppose that $p=2$. For $n\geq 1$, we denote by $G_n$ the (almost) extraspecial 2-group described in Lemma 2.1 \textit{i)}, \textit{ii)}. Hence, $G_1\cong D_8$ or $G_1\cong Q_8$ in the extraspecial case, and $G_1\cong D_8\circ C_4$ in the almost extraspecial case. Also, for $n\geq 2$, each (almost) extraspecial $2$-group $G_{n}$ can be written as $D_8\circ G_{n-1}$. Under this notation, the author proved the following result in \cite{22}.\\

\textbf{Lemma 2.4.} \textit{Let $G_{n}\cong D_8\circ G_{n-1}$ be an (almost) extraspecial 2-group of order $2^{m}$ where $m\in \lbrace 2n+1, 2n+2 \rbrace$ and $n\geq 2$ are integers. Then $n_2(G_{n})=2^{m-2}+2n_2(G_{n-1})+1.$}\\

Let $G$ be a (generalized/almost) extraspecial $2$-group. Since $|G|$ is given by Lemma 2.1 and $\exp(G)=4$ by Lemma 2.2 \textit{i)}, it follows that once we determine $n_2(G)$, we can easily find $n_4(G)$ and $\psi(G)$. The following theorem highlights such numerical results. For the ease of writing, we use the same notation as in Lemma 2.1.\\  

\textbf{Theorem 2.5.} \textit{Let $G$ be a (generalized/almost) extraspecial 2-group.
\begin{itemize}
\item[i)] If $G\cong D_8^{\circ n}$, then $n_2(G)=4^n+2^n-1, n_4(G)=4^n-2^n$ and $$\psi(G)=6\cdot 2^{2n}-2\cdot 2^n-1;$$
\item[ii)] If $G\cong D_8^{\circ (n-1)}\circ Q_8$, then $n_2(G)=4^n-2^n-1, n_4(G)=4^n+2^n$ and $$\psi(G)=6\cdot 2^{2n}+2\cdot 2^n-1;$$
\item[iii)] If $G\cong D_8^{\circ n}\circ C_4$, then $n_2(G)=2^{2n+1}-1, n_4(G)=2^{2n+1}$ and $$\psi(G)=12\cdot 2^{2n}-1;$$
\item[iv)] If $G\cong D_8^{\circ n}\times C_2^k$, then $n_2(G)=2^{2n+k}+2^{n+k}-1, n_4(G)=2^{2n+k}-2^{n+k}$ and $$\psi(G)=6\cdot 2^{2n+k}-2\cdot 2^{n+k}-1;$$
\item[v)] If $G\cong (D_8^{\circ (n-1)}\circ Q_8)\times C_2^k$, then $n_2(G)=2^{2n+k}-2^{n+k}-1, n_4(G)=2^{2n+k}+2^{n+k}$ and $$\psi(G)=6\cdot 2^{2n+k}+2\cdot 2^{n+k}-1;$$
\item[vi)] If $G\cong (D_8^{\circ n}\circ C_4)\times C_2^k$, then $n_2(G)=2^{2n+k+1}-1, n_4(G)=2^{2n+k+1}$ and $$\psi(G)=12\cdot 2^{2n+k}-1.$$
\end{itemize}}

\textbf{Proof.} Using the previous established notations, let $G\cong G_n$ be an (almost) extraspecial $2$-group of order $2^m$ where $m\in \lbrace 2n+1, 2n+2\rbrace$ and $n\geq 1$ are integers. If $G_1$ is extraspecial, then $n_2(G_1)=n_2(D_8)=5$ or $n_2(G_1)=n_2(Q_8)=1$. If $G_1$ is almost extraspecial, then $n_2(G_1)=n_2(D_8\circ C_4)=7$. The same results are obtained by using the formulas outlined  by items \textit{i)}, \textit{ii)} and \textit{iii)} for $n=1$.

We show by induction on $n$ that
\begin{align}\label{r1}
n_2(G_n)=2^{m-n}(2^{n-1}-1)+2^{n-1}n_2(G_1)+2^{n-1}-1, \ \forall \ n\geq 2.
\end{align}
The base case is exactly Lemma 2.4 applied for $n=2$. Let $n\geq 3$ and suppose that (\ref{r1}) holds for any integer $k$ such that $2\leq k<n$. Then, by applying Lemma 2.4 and the inductive hypothesis, we have
\begin{align*}
n_2(G_{n})&=2^{m-2}+2n_2(G_{n-1})+1\\ &=2^{m-2}+2[2^{m-n-1}(2^{n-2}-1)+2^{n-2}n_2(G_1)+2^{n-2}-1]+1\\ &=2^{m-n}(2^{n-1}-1)+2^{n-1}n_2(G_1)+2^{n-1}-1,
\end{align*}
as stated.

If $m=2n+1$, then $G_n$ is extraspecial and relation (\ref{r1}) becomes
\begin{align}\label{r2}
n_2(G_n)=2^{n+1}(2^{n-1}-1)+2^{n-1}n_2(G_{1})+2^{n-1}-1.
\end{align}
In this case, $G_1\cong D_8$ or $G_1\cong Q_8$. Since $n_2(D_8)=5$ and $n_2(Q_8)=1$, by making the replacements in (\ref{r2}), we get the quantities $n_2(G)$ outlined by items \textit{i)} and \textit{ii)} of this theorem.

If $m=2n+2$, then $G_n$ is almost extraspecial. In this case, relation (\ref{r1}) may be rewritten as
\begin{align}\label{r3}
n_2(G_n)=2^{n+2}(2^{n-1}-1)+2^{n-1}n_2(G_1)+2^{n-1}-1.
\end{align}
Again, since $G_1\cong D_8\circ C_4$ and $n_2(D_8\circ C_4)=7$, relation (\ref{r3}) leads to obtaining the value of $n_2(G)$ given by item \textit{iii)}.

Assume now that $G$ is a generalized extraspecial $2$-group. Then, for items \textit{iv)}, \textit{v)} and \textit{vi)}, the number $n_2(G)$ is determined based on the number of elements of order $2$ of an (almost) extraspecial group and by applying Proposition 2.3. For instance, if we are to refer to item \textit{iv)}, we have
\begin{align*}
n_2(G)=n_2(D_8^{\circ n}\times C_2^k)=  (n_2(D_8^{\circ n})+1)(n_2(C_2^k)+1)-1=(4^n+2^n)2^k-1=2^{2n+k}+2^{n+k}-1,
\end{align*}
as desired. The same reasoning is applied in order to obtain the quantity $n_2(G)$ highlighted by items \textit{v)} and \textit{vi)}.

Finally, for all $6$ main cases, the number of elements of order $4$ and the sum of element orders of $G$ are given by
$$n_4(G)=|G|-1-n_2(G) \text{ \ \ \ \ \ \ \ \ \ \ and \ \ \ \ \ \ \ \ \ \ } \psi(G)=1+2n_2(G)+4n_4(G),$$
respectively. Thus, the proof is complete.
\hfill\rule{1,5mm}{1,5mm}\\    

In the second part of this section, we are mainly interested in investigating the element orders of a (generalized/almost) extraspecial $p$-group where $p$ is odd. Still, we start by justifying a simple, yet useful, result which is also valid for $p=2$.\\

\textbf{Lemma 2.6.} \textit{Let $G$ be a $p$-group of order $p^n$ where $n\geq 3$ is an integer. Let $k\leq n$ and $j\leq n$ be non-negative integers. If $\Omega_k(G)=\Omega_{\lbrace k\rbrace}(G)$ and $[G:\Omega_k(G)]=p^j$, then $|\Omega_{\lbrace k\rbrace}(G)|=p^{n-j}$.}\\

\textbf{Proof.} Under the highlighted hypotheses, the condition $[G:\Omega_k(G)]=p^j$ implies that $|\Omega_k(G)|=p^{n-j}$. Since $\Omega_k(G)=\Omega_{\lbrace k\rbrace}(G)$, the conclusion follows.
\hfill\rule{1,5mm}{1,5mm}\\    

We recall that a finite $p$-group $G$ is called regular if, given any positive integer $n$ and any $a, b\in G$, there are $c_3, c_4, \ldots, c_k\in \langle a, b\rangle'$ such that $(ab)^{p^n}=a^{p^n}b^{p^n}c_3^{p^n}c_4^{p^n}\ldots c_k^{p^n}$. Also, a finite $p$-group $G$ is called powerful if either $p$ is odd and $G'\leq \mho^1(G)$ or $p=2$ and $G'\leq \mho^2(G)$.
For a regular $p$-group $G$, it is known that $\Omega_k(G)=\Omega_{\lbrace k\rbrace}(G)$, for all $k\geq 1$ (see section 4.3 of \cite{13}). Hence, for such groups, if one knows the values of the indices $[G:\Omega_k(G)]$ (or, equivalently, $|\mho^k(G)|$; see Theorem 7.2 \textit{(d}) of \cite{6}), then Lemma 2.6 leads to an easy way to determine the numbers $n_{p^k}(G)$ based on the formula
\begin{align*}
n_{p^k}(G)=|\Omega_{\lbrace k\rbrace}(G)|-|\Omega_{\lbrace k-1\rbrace}(G)|, \ \forall \ k\geq 1.
\end{align*} 
The same thing can be said about powerful $p$-groups where $p$ is odd (for a powerful $2$-group, the equality $\Omega_k(G)=\Omega_{\lbrace k\rbrace}(G)$ does not hold in general; indeed, if $G\cong D_8\circ C_4$, we have $|\Omega_{\lbrace 1\rbrace}(G)|=8$, but $|\Omega_{1}(G)|=16$). For more details on such properties of powerful $p$-groups, one may check \cite{16, 23, 24}. 

A (generalized/almost) extraspecial $p$-group $G$ where $p$ is odd, is always regular. This follows from the fact that $G$ is of odd order and $G'$ is cyclic (see 10.2 Satz \textit{c)} of \cite{17}). Also, except the case in which $\exp(G)=p$, we have that $G$ is powerful by the definition of the powerfulness property. Indeed, since $\lbrace 1\rbrace\ne\mho^1(G)\leq \Phi(G)$, $\Phi(G)=G'$ and $|G'|=p$, we deduce that $G'=\mho^1(G)$, so $G$ is powerful. 

The following result completes our study on the element orders of a (generalized/almost) extraspecial $p$-group. Again, for the ease of writing, we use the same notation as in Lemma 2.1.\\  

\textbf{Theorem 2.7.} \textit{Let $G$ be a (generalized/almost) extraspecial $p$-group where $p$ is odd.
\begin{itemize}   
\item[i)] If $G\cong He_p^{\circ n}$, then $n_p(G)=p^{2n+1}-1$ and $\psi(G)=p^{2n+2}-p+1$;
\item[ii)] If $G\cong M_{p^3}^{\circ n}$, then $n_p(G)=p^{2n}-1, n_{p^2}(G)=(p-1)p^{2n}$ and $$\psi(G)=(p-1)p^{2n+2}+p(p^{2n}-1)+1;$$
\item[iii)] If $G\cong He_p^{\circ n}\circ C_{p^2}$, then $n_p(G)=p^{2n+1}-1, n_{p^2}(G)=(p-1)p^{2n+1}$ and $$\psi(G)=(p-1)p^{2n+3}+p(p^{2n+1}-1)+1;$$
\item[iv)] If $G\cong He_p^{\circ n}\times C_p^k$, then $n_p(G)=p^{2n+k+1}-1$ and $\psi(G)=p^{2n+k+2}-p+1$;
\item[v)] If $G\cong M_{p^3}^{\circ n}\times C_p^k$, then $n_p(G)=p^{2n+k}-1, n_{p^2}(G)=(p-1)p^{2n+k}$ and $$\psi(G)=(p-1)p^{2n+k+2}+p(p^{2n+k}-1)+1;$$
\item[vi)] If $G\cong (He_p^{\circ n}\circ C_{p^2})\times C_p^k$, then $n_p(G)=p^{2n+k+1}-1, n_{p^2}(G)=(p-1)p^{2n+k+1}$ and $$\psi(G)=(p-1)p^{2n+k+3}+p(p^{2n+k+1}-1)+1.$$
\end{itemize}}

\textbf{Proof.} Let $G$ be a (generalized/almost) extraspecial $p$-group where $p$ is odd. For items \textit{i)} and \textit{iv)}, by Lemma 2.2 \textit{ii)} we have $\exp(G)=p$ and, by using the information on $|G|$ given by Lemma 2.1, we easily determine $n_p(G)$ and $\psi(G)$.

For all the other 4 cases, we know that $\exp(G)=p^2$ by Lemma 2.2 \textit{ii)}. Since $G$ is regular, we have $\Omega_1(G)=\Omega_{\lbrace 1\rbrace}(G)$ and $[G:\Omega_1(G)]=|\mho^1(G)|$. As we previously noticed, $\mho^1(G)=G'$ and, since $|G'|=p$, by applying Lemma 2.6 for $k=j=1$, we get $|\Omega_{\lbrace 1\rbrace}(G)|=\frac{|G|}{p}$. Hence $n_p(G)=\frac{|G|}{p}-1$. By replacing $|G|$ with $p^{2n+1}$, $p^{2n+2}$, $p^{2n+k+1}$ and $p^{2n+k+2}$, we obtain the values of $n_p(G)$ given by items \textit{ii)}, \textit{iii)}, \textit{v)} and \textit{vi)}, respectively. To conclude, for the same 4 cases, we can easily deduce the number of elements of order $p^2$ and the sum of element orders by applying the formulas
$$n_{p^2}(G)=|G|-1-n_p(G) \text{ \ \ \ \ \ \ \ \ \ \ and \ \ \ \ \ \ \ \ \ \ } \psi(G)=1+pn_p(G)+p^2n_{p^2}(G),$$
respectively.
\hfill\rule{1,5mm}{1,5mm}\\ 

For the last 3 cases of Theorem 2.7, another approach to determine $n_p(G)$ is to use Proposition 2.3, as we did in the proof of Theorem 2.5. Some examples of applying the two main theorems are given below. We mention that these results are also supported by GAP \cite{35}.\\

\textbf{Example 2.8.} \textit{We have}\\

\begin{tabular}{ |p{3cm}|p{1cm}|p{2cm}|p{2cm}|p{2cm}|p{2cm}|  }
 \hline
  Group & Order & In GAP & $n_p(G)$ & $n_{p^2}(G)$ & $\psi(G)$\\
 \hline
 $D_8^{\circ 3}$   &  128  & 2326 & 71 & 56 & 367 \\
 \hline
 $D_8^{\circ 2}\circ Q_8$ & 128  & 2327   & 55 & 72 & 399 \\
 \hline
  $D_8^{\circ 2}\circ C_4$ & 64  & 266   & 31 & 32 & 191 \\
 \hline
 $ D_8^{\circ 2}\times C_2^2$ & 128  & 2323   & 79 & 48 & 351 \\
 \hline
 $ (D_8\circ Q_8)\times C_2^2$ & 128  & 2324   & 47 & 80 & 415 \\
 \hline
 $ (D_8^{\circ 2}\circ C_4)\times C_2^2$ & 256  & 56088   & 127 & 128 & 767 \\
 \hline
 $ He_5^{\circ 3}$ & 78125 & 34295   & 78124 & - & 390621 \\
 \hline
 $ M_{125}^{\circ 3}$ & 78125 & 34296   & 15624 & 62500 & 1640621 \\
 \hline
 $ He_5^{\circ 2}\circ C_{25}$ & 15625 & 122   & 3124 & 12500 & 328121 \\
 \hline
 $ He_5^{\circ 2}\times C_5^2$ & 78125 & 34292   & 78124 & - & 390621 \\
 \hline
  $ M_{125}^{\circ 2}\times C_5^2$ & 78125 & 34293   & 15624 & 62500 & 1640621\\
 \hline
 $ (He_5\circ C_{25})\times C_5^2$ & 15625 & 26   & 3124 & 12500 & 328121\\
 \hline
\end{tabular}
\vspace{0.2in}

The last 6 lines of the previous table show that a (generalized/almost) extraspecial $p$-group $G$ where $p$ is odd, is not uniquely determined by $|G|$ and $\psi(G)$. The same can be said for $p=2$ since  
$$\begin{array}{ll}
\psi(D_8^{\circ 3})=\psi(D_8^2\times C_2)=367; \\
\psi(D_8\circ Q_8)=\psi(D_8\times C_4)=103; \\
\psi(D_8\circ C_4)=\psi(C_2^2\times C_4)=47;
\end{array} \ \ \  
\begin{array}{ll}
\psi(D_8^{\circ 3}\times C_2)=\psi(D_8^2\times C_2^2)=735; \\
\psi(Q_8\times C_2)=\psi(C_4^2)=55; \\
\psi((D_8\circ C_4)\times C_2)=\psi(C_2^3\times C_4)=95. 
\end{array}$$
  
We end this section by pointing out that, as consequences of Theorems 2.5 and 2.7, one can easily obtain the number of cyclic subgroups of all possible orders for any (generalized/almost) extraspecial $p$-group.  

\section{A density result on the cyclicity degree of a finite group}

Our following main aim is to give a positive answer to Question 1.1 which was recalled in the first section of this paper. Along with the notations already established, for a positive integer $k\geq 1$, we denote by $p_k$ the \textit{k}th odd prime number. Also, for a topological space $(X, \tau_X)$ and $A\subseteq X$, the closure of $A$ with respect to $\tau_X$ is denoted by $\overline{A}_{\tau_X}$. We start by highlighting some preliminary results related to calculus and to the cyclicity degree of a finite group.\\ 

\textbf{Lemma 3.1.} \textit{Let $(x_k)_{k\geq 1}, (y_k)_{k\geq 1}$ be sequences of positive real numbers. Let $(G_k)_{k\in I}$ be a family of finite groups where $I$ is a finite non-empty set. Let $(X, \tau_X)$ and $(Y, \tau_Y)$ be topological spaces, $f\colon (X, \tau_X)\longrightarrow (Y, \tau_Y)$ be a continuous function and $A, B\subseteq X$.
\begin{itemize}
\item[i)] If $\displaystyle\lim_{k\to\infty}\frac{x_k}{y_k}\in (0, \infty),$
then the series $\sum\limits_{k=1}^{\infty}x_k$ and $\sum\limits_{k=1}^{\infty}y_k$ have the same nature (see Theorem 10.9 of \cite{2});
\item[ii)] $\sum\limits_{k=1}^{\infty}\frac{1}{p_k}=\infty$ (consequence of the main result of \cite{25});
\item[iii)] If $\displaystyle\lim_{k\to\infty}x_k=0$ and $\sum\limits_{k=1}^{\infty}x_k=\infty$,
then the set containing the sums of all finite subsequences of $(x_k)_{k\geq 1}$ is dense in $[0, \infty)$ (consequence of the Proposition outlined on p. 863 of \cite{26});
\item[iv)] If $\overline{A}_{\tau_X}=\overline{B}_{\tau_X}$, then $\overline{f(A)}_{\tau_Y}=\overline{f(B)}_{\tau_Y}$ (see Proposition 6.12 of \cite{31});
\item[v)] $\cdeg(M_{p^n})=\frac{(n-1)p+2}{(n-1)p+n+1}$ (see Theorem 3.3.2 of \cite{34});
\item[vi)] $\cdeg\bigg(\xmare{k\in I}{ }G_{k}\bigg)=\prod\limits_{k\in I}\cdeg(G_k)$ (see Proposition 2.2 of \cite{34}).
\end{itemize}}
\vspace{0.2in}

For a class $\mathcal{C}$ of finite groups, we denote by 
$$CD_{\mathcal{C}}=\lbrace \cdeg(G) \ | \ G\in\mathcal{C}\rbrace.$$

The following theorem is the main result of this section and it guarantees that Question 1.1 can be answered positively.\\

\textbf{Theorem 3.2.} \textit{Let $\mathcal{M}$ be a class of finite groups containing all finite direct products of distinct  extraspecial $p$-groups of the form $M_{p^3}$ where $p$ is an odd prime. Then $CD_{\mathcal{M}}$ is dense in $[0, 1].$}\\

\textbf{Proof.} Let $(y_k)_{k\geq 1}$ be the sequence of the reciprocals of the odd primes, i.e.
$$y_k=\frac{1}{p_k}, \ \forall \ k\geq 1.$$
We define a sequence of positive real numbers $(x_k)_{k\geq 1}$ where
$$x_k=\ln\frac{1}{\cdeg(M_{p_k^3})}, \ \forall \ k\geq 1.$$
By Lemma 3.1 \textit{v)} applied for $n=3$, we deduce that
\begin{align}\label{r5}
\displaystyle\lim_{k\to\infty}x_k=\displaystyle\lim_{k\to\infty}\ln\bigg(\frac{p_k+2}{p_k+1}\bigg)=0.
\end{align}
Also, since
$$\displaystyle\lim_{k\to\infty}\frac{x_k}{y_k}=\displaystyle\lim_{k\to\infty}\frac{\ln\big(\frac{p_k+2}{p_k+1}\big)}{\frac{1}{p_k}}=1,$$
by Lemma 3.1 \textit{i), ii)}, we have
\begin{align}\label{r6}
\sum\limits_{k=1}^{\infty}x_k=\infty.
\end{align}
Hence, by (\ref{r5}) and (\ref{r6}), we conclude that the sequence $(x_k)_{k\geq 1}$ satisfies the hypotheses of Lemma 3.1 \textit{iii)}. This implies that
$$\overline{\bigg\lbrace \sum\limits_{k\in I}x_k \ \bigg| \ I\subseteq \mathbb{N}^*, |I|<\infty\bigg\rbrace}_{\tau_{\mathbb{R}}}=[0, \infty),$$
so, by the standard properties of the natural logarithm and Lemma 3.1 \textit{vi)}, we have
\begin{align}\label{r7}
\overline{\Bigg\lbrace \ln\frac{1}{\cdeg\bigg(\xmare{k\in I}{ } M_{p_k^3}\bigg)} \ \Bigg| \ I\subseteq \mathbb{N}^*, |I|<\infty\Bigg\rbrace}_{\tau_{\mathbb{R}}}=[0, \infty),
\end{align}
where $\tau_{\mathbb{R}}$ is the usual topology of $\mathbb{R}$. Since (\ref{r7}) outlines an equality between the closures of two subsets of $\mathbb{R}$ and the exponential function
$$\exp:(\mathbb{R}, \tau_{\mathbb{R}})\longrightarrow (\mathbb{R}, \tau_{\mathbb{R}}), \text{ given by } \exp(x)=e^x, \ \forall \ x\in \mathbb{R},$$
is continuous, by Lemma 3.1 \textit{iv)}, it follows that
\begin{align}\label{r8}
\overline{\Bigg\lbrace \frac{1}{\cdeg\bigg(\xmare{k\in I}{ } M_{p_k^3}\bigg)} \ \Bigg| \ I\subseteq \mathbb{N}^*, |I|<\infty\Bigg\rbrace}_{\tau_{\mathbb{R}}}=[1, \infty).
\end{align}

Let $Y=(0, \infty)$ and $\tau_Y$ be the subspace topology on $Y$. It is known that for a subset $X$ of $Y$, we have $\overline{X}_{\tau_{Y}}=\overline{X}_{\tau_{\mathbb{R}}}\cap Y$. Hence, by (\ref{r8}), we deduce that
\begin{align}\label{r9}
\overline{\Bigg\lbrace \frac{1}{\cdeg\bigg(\xmare{k\in I}{ } M_{p_k^3}\bigg)} \ \Bigg| \ I\subseteq \mathbb{N}^*, |I|<\infty\Bigg\rbrace}_{\tau_Y}=[1, \infty).
\end{align}
Finally, by applying the continuous function 
$$f:(Y, \tau_Y)\longrightarrow (\mathbb{R}, \tau_{\mathbb{R}}), \text{ given by } f(y)=\frac{1}{y}, \ \forall \ y\in Y,$$
to both sides of (\ref{r9}), we get
\begin{align}\label{r10}
\overline{\bigg\lbrace \cdeg\bigg(\xmare{k\in I}{ } M_{p_k^3}\bigg) \bigg| \ I\subseteq \mathbb{N}^*, |I|<\infty\bigg\rbrace}_{\tau_{\mathbb{R}}}=[0, 1].
\end{align}
Since 
$$\bigg\lbrace \cdeg\bigg(\xmare{k\in I}{ } M_{p_k^3}\bigg) \bigg| \ I\subseteq \mathbb{N}^*, |I|<\infty\bigg\rbrace\subseteq CD_{\mathcal{M}}\subseteq [0, 1],$$
and taking closures preserve inclusion, we get
$$\overline{\bigg\lbrace \cdeg\bigg(\xmare{k\in I}{ } M_{p_k^3}\bigg) \bigg| \ I\subseteq \mathbb{N}^*, |I|<\infty\bigg\rbrace}_{\tau_{\mathbb{R}}}\subseteq \overline{CD_{\mathcal{M}}}_{\tau_{\mathbb{R}}} \subseteq [0, 1].$$
Therefore, as a consequence of (\ref{r10}), we obtain
$$\overline{CD_{\mathcal{M}}}_{\tau_{\mathbb{R}}}=[0, 1],$$
so,  $CD_{\mathcal{M}}$ is dense in $[0, 1]$.
\hfill\rule{1,5mm}{1,5mm}\\ 

Theorem 3.2 implies the following immediate result which end our paper. Its second item constitutes an affirmative answer to Question 1.1.\\

\textbf{Corollary 3.3.} \textit{\begin{itemize}
\item[i)] Let $\mathcal{M}$ be a class of finite groups containing all finite direct products of distinct extraspecial $p$-groups of the form $M_{p^3}$ where $p$ is an odd prime, and let $\mathcal{C}$ be a class of finite groups such that $\mathcal{M}\subseteq\mathcal{C}$. Then $CD_{\mathcal{C}}$ is dense in $[0, 1].$
\item[ii)] Let $a\in [0, 1]$. Then there exists a sequence $(G_n)_{n\geq 1}$ of finite groups such that $$\displaystyle\lim_{n\to\infty} \cdeg(G_n)=a.$$ 
\end{itemize}}

Based on the proof of Theorem 3.2, each term of the sequence $(G_n)_{n\geq 1}$ is a finite direct product of distinct extraspecial $p$-groups of the form $M_{p^3}$ where $p$ is an odd prime. \\

\bigskip\noindent {\bf Acknowledgements.} The author is grateful to the reviewer for his/her remarks which improved the initial version of the paper.

\vspace*{3ex}
\small
\hfill
\begin{minipage}[t]{7cm}
Mihai-Silviu Lazorec \\
Faculty of  Mathematics \\
"Al.I. Cuza" University \\
Ia\c si, Romania \\
e-mail: {\tt silviu.lazorec@uaic.ro}
\end{minipage}
\end{document}